\documentclass[12pt]{article}
\usepackage{amsfonts}
\usepackage{a4wide}
\usepackage{amsmath,amscd,amsthm,a4,amssymb}

\begin{document}

\begin{center}
{\Large \textbf{A note on the equioscillation theorem for best ridge
function approximation}}

\

\textbf{Vugar E. Ismailov} \

\smallskip

Institute of Mathematics and Mechanics

National Academy of Sciences of Azerbaijan

Az-1141, Baku, Azerbaijan

{e-mail:} {vugaris@mail.ru}
\end{center}

\smallskip

\textbf{Abstract.} We consider the approximation of a continuous function,
defined on a compact set of the $d$-dimensional Euclidean space, by sums of
two ridge functions. We obtain a necessary and sufficient condition for such a
sum to be a best approximation. The result resembles the classical Chebyshev
equioscillation theorem for polynomial approximation.

\bigskip



\bigskip

\begin{center}
{\large \textbf{1. Introduction}}
\end{center}

Let $Q$ be compact set in the $d$-dimensional Euclidean space and $C(Q)$ be
the space of continuous real-valued functions on $Q$. Consider the
approximation of a function $f\in C(Q)$ by sums of the form $g_{1}(\mathbf{a}%
_{1}\cdot \mathbf{x})+g_{2}(\mathbf{a}_{2}\cdot \mathbf{x}),$ where $\mathbf{%
a}_{i}$ are fixed vectors (directions) in $\mathbb{R}^{d}\backslash \left\{
\mathbf{0}\right\} $ and $g_{i}$ are continuous univariate functions. We are
interested in characterization of a best approximation. Note that functions
of the form $g(\mathbf{a}\cdot \mathbf{x})$ are called ridge functions.
These functions and their linear combinations arise naturally in problems of
computerized tomography (see, e.g., \cite{Log}), statistics (see, e.g., \cite%
{Fri,Hub}), partial differential equations \cite{John} (where they are
called \textit{plane waves}), neural networks (see, e.g., \cite{P2} and
references therein), and approximation theory (see, e.g., \cite{7,17,Pel,Pet}%
). In the past few years, problems of ridge function representation have
gained special attention among researchers (see e.g. \cite{Ali,Kon1,Kon2,Pin}%
). For more on ridge functions and application areas see a recently
published monograph by Pinkus \cite{P1}.

Characterization theorems for best approximating elements are essential in
approximation theory. The classical and most striking example of such a
theorem are the Chebyshev equioscillation theorem. This theorem
characterizes the unique best uniform approximation to a continuous real
valued function $F(t)$ by polynomials $P(t)$ of degree at most $n$, by the
oscillating nature of the difference $F(t)-P(t)$. The result says that if
such polynomial has the property that for some particular $n+2$ points $%
t_{i} $ in $[0,1]$

\begin{equation*}
F(t_{i})-P(t_{i})=(-1)^{i}\max_{x\in \lbrack 0,1]}\left\vert
F(t)-P(t)\right\vert ,\text{ }i=1,...,n+2,
\end{equation*}%
then $P$ is the best approximation to $F$ on $[0,1]$. The monograph of
Natanson \cite{Nat} contains a very rich commentary on this theorem. Some
general alternation type theorems applying to any finite dimensional
subspace $E$ of $C(I)$ for $I$ a cell in $\mathbb{R}^{d}$, may be found in
Buck \cite{Buck}. For a short history and various modifications of the
Chebyshev alternation theorem see \cite{Bro}.

In this note, we obtain an equioscillation theorem for approximation of
multivariate functions by sums of two ridge functions. To be more precise,
let $Q$ be a compact subset of the space $\mathbb{R}^{d}$. Fix two
directions $\mathbf{a}_{1}$ and $\mathbf{a}_{2}$ in $\mathbb{R}^{d}$ and
consider the following space
\begin{equation*}
\mathcal{R}=\mathcal{R}(\mathbf{a}_{1},\mathbf{a}_{2})=\{g_{1}(\mathbf{a}%
_{1}\cdot \mathbf{x})+g_{2}(\mathbf{a}_{2}\cdot \mathbf{x}):~g_{1},g_{2}\in
C(\mathbb{R})\}.
\end{equation*}%
Note that the space $\mathcal{R}$ is a linear space. Assume a function $f\in
C(Q)$ is given. We ask and answer the following question: which geometrical
conditions imposed on $G_{0}\in \mathcal{R}$ is necessary and sufficient for
the equality

\begin{equation*}
\left\Vert f-G_{{0}}\right\Vert =\inf_{G\in \mathcal{R}}\left\Vert
f-G\right\Vert \text{?}\eqno(1.1)
\end{equation*}%
Here $\left\Vert \cdot \right\Vert $ denotes the standard uniform norm in $%
C(Q).$ Recall that functions $G_{0}$ satisfying (1.1) are called best
approximations or extremal elements.

It should be remarked that in the special case when $Q\subset \mathbb{R}^{2}$
and $\mathbf{a}_{1}$ and $\mathbf{a}_{2}$ coincide with the coordinate
directions, the above question was answered by Khavinson \cite{Kh}. In \cite%
{Kh}, he obtained an equioscillation theorem for a best approximating sum $%
\varphi (x)+\psi (y)$. In our papers \cite{7,11}, Chebyshev type theorems
were proven for ridge functions under additional assumption that $Q$ is
convex. For a more recent and detailed discussion of an equioscillation
theorem in ridge function approximation see Pinkus \cite{P1}.

\bigskip

\bigskip

\begin{center}
{\large \textbf{2. Equioscillation theorem for ridge functions}}
\end{center}

We start with a definition of paths with respect to two directions. These
objects will play an essential role in our further analysis.

\bigskip

\textbf{Definition 2.1 (see \cite{11}).} \ \textit{A finite or infinite
ordered set $p=\left( \mathbf{p}{_{1},\mathbf{p}_{2},...}\right) \subset Q$
with $\mathbf{p}_{i}\neq \mathbf{p}_{i+1},$ and either $\mathbf{a}_{1}\cdot
\mathbf{p}_{1}=\mathbf{a}_{1}\cdot \mathbf{p}_{2},\mathbf{a}_{2}\cdot
\mathbf{p}_{2}=\mathbf{a}_{2}\cdot \mathbf{p}_{3},\mathbf{a}_{1}\cdot
\mathbf{p}_{3}=\mathbf{a}_{1}\cdot \mathbf{p}_{4},...$ or $\mathbf{a}%
_{2}\cdot \mathbf{p}_{1}=\mathbf{a}_{2}\cdot \mathbf{p}_{2},~\mathbf{a}%
_{1}\cdot \mathbf{p}_{2}=\mathbf{a}_{1}\cdot \mathbf{p}_{3},\mathbf{a}%
_{2}\cdot \mathbf{p}_{3}=\mathbf{a}_{2}\cdot \mathbf{p}_{4},...$is called a
path with respect to the directions $\mathbf{a}_{1}$ and $\mathbf{a}_{2}$.}

\bigskip

In the sequel, we will simply use the term \textquotedblleft
path\textquotedblright\ instead of the expression \textquotedblleft path
with respect to the directions $\mathbf{a}_{1}$ and $\mathbf{a}_{2}$%
\textquotedblright . If in a finite path $(\mathbf{p}_{1},...,\mathbf{p}_{n},%
\mathbf{p}_{n+1})$, $\mathbf{p}_{n+1}=\mathbf{p}_{1}$ and $n$ is an even
number, then the path $(\mathbf{p}_{1},...,\mathbf{p}_{n})$ is said to be
closed. Note that for a closed path $(\mathbf{p}_{1},...,\mathbf{p}_{2n})$
and any function $G\in \mathcal{R}$, $G(\mathbf{p}_{1})-G(\mathbf{p}%
_{2})+\cdot \cdot \cdot -G(\mathbf{p}_{2n})=0$.

Paths, in the special case when $Q\subset \mathbb{R}^{2}$, $\mathbf{a}_{1}$
and $\mathbf{a}_{2}$ coincide with the coordinate directions, are
geometrically explicit objects. In this case, a path is a finite ordered set
$(\mathbf{p}_{1},...,\mathbf{p}_{n})$ in $\mathbb{R}^{2}$ with the line
segments $[\mathbf{p}_{i},\mathbf{p}_{i+1}],$ $i=1,...,n,$ alternatively
perpendicular to the $x$ and $y$ axes (see, e.g., \cite%
{Arn,5,6,Is2,Is3,12,Mar1}). These objects were first introduced by Diliberto
and Straus \cite{4} (in \cite{4}, they are called \textquotedblleft
permissible lines\textquotedblright ). They appeared further in a number of
papers with several different names such as \textquotedblleft bolts" (see,
e.g., \cite{Arn,12,Mar1}), \textquotedblleft trips\textquotedblright\ (see
\cite{Mar2}), \textquotedblleft links" (see, e.g., \cite{Cow,Klo1,Klo2}),
etc. Paths with respect to two directions $\mathbf{a}_{1}$ and $\mathbf{a}%
_{2}$ were exploited in some papers devoted to ridge function interpolation
(see, e.g., \cite{BP,IP}). In \cite{9,10}, paths were generalized to those
with respect to a finite set of functions. The last objects turned out to be
very useful in problems of representation by linear superpositions.

In the sequel, we need the concept of an \textquotedblleft extremal path",
which is defined as follows.

\bigskip

\textbf{Definition 2.2 (see \cite{11}).}\ \textit{A finite or infinite path $%
(\mathbf{p}_{1},\mathbf{p}_{2},...)$ is said to be extremal for a function $%
h\in C(Q)$ if $h(\mathbf{p}_{i})=(-1)^{i}\left\Vert h\right\Vert ,i=1,2,...$
or $h(\mathbf{p}_{i})=(-1)^{i+1}\left\Vert h\right\Vert ,$ $i=1,2,...$}

\bigskip

The purpose of this note is to prove the following theorem.

\bigskip

\textbf{Theorem 2.1.}\ \textit{\ Assume $Q$ is a compact subset of $\mathbb{R%
}^{d}$. A function $G_{0}\in \mathcal{R}$ is a best approximation to a
function $f\in C(Q)$ if and only if there exists a closed or infinite path $%
p=(\mathbf{p}_{1},\mathbf{p}_{2},...)$ extremal for the function $f-G_{0}$.}

\bigskip

\textbf{Proof. Sufficiency. }There are two possible cases. The first case
happens when there exists a closed path $(\mathbf{p}_{1},...,\mathbf{p}%
_{2n}) $ extremal for the function $f-G_{0}.$ Let us check that in this
case, $f-G_{0}$ is a best approximation. Indeed, on the one hand, the
following equalities are valid.

\begin{equation*}
\left\vert \sum_{i=1}^{2n}(-1)^{i}f(\mathbf{p}_{i})\right\vert =\left\vert
\sum_{i=1}^{2n}(-1)^{i}\left[ f-G_{0}\right] (\mathbf{p}_{i})\right\vert
=2n\left\Vert f-G_{0}\right\Vert .
\end{equation*}%
On the other hand, for any function $G\in \mathcal{R}$, we have

\begin{equation*}
\left\vert \sum_{i=1}^{2n}(-1)^{i}f(\mathbf{p}_{i})\right\vert =\left\vert
\sum_{i=1}^{2n}(-1)^{i}\left[ f-G\right] (\mathbf{p}_{i})\right\vert \leq
2n\left\Vert f-G\right\Vert .
\end{equation*}%
Therefore, $\left\Vert f-G_{0}\right\Vert \leq \left\Vert f-G\right\Vert $
for any $G\in \mathcal{R}$. That is, $G_{0}$ is a best approximation.

The second case happens when we do not have closed paths extremal for $%
f-G_{0}$, but there exists an infinite path $(\mathbf{p}_{1},\mathbf{p}%
_{2},...)$ extremal for $f-G_{0}$. To analyze this case, consider the
following linear functional

\begin{equation*}
l_{q}:C(Q)\rightarrow \mathbb{R}\text{, \ }l_{q}(F)=\frac{1}{n}%
\sum_{i=1}^{n}(-1)^{i}F(\mathbf{q}_{i}),
\end{equation*}%
where $q=\{\mathbf{q}_{1},...,\mathbf{q}_{n}\}$ is a finite path in $Q$. It
is easy to see that the norm $\left\Vert l_{q}\right\Vert \leq 1$ and $%
\left\Vert l_{q}\right\Vert =1$ if and only if the set of points of $q$ with
odd indices $O=\{\mathbf{q}_{i}\in q:$ $i$ \textit{is an odd number}$\}$ do
not intersect with the set of points of $q$ with even indices $E=\{\mathbf{q}%
_{i}\in q:$ $i$ \textit{is an even number}$\}$. Indeed, from the definition
of $l_{q}$ it follows that $\left\vert l_{q}(F)\right\vert \leq \left\Vert
F\right\Vert $ for all functions $F\in C(Q)$, whence $\left\Vert
l_{q}\right\Vert \leq 1.$ If $O\cap E=\varnothing $, then for a function $%
F_{0}$ with the property $F_{0}(\mathbf{q}_{i})=-1$ if $i$ is odd, $F_{0}(%
\mathbf{q}_{i})=1$ if $i$ is even and $-1<F_{0}(x)<1$ elsewhere on $Q,$ we
have $\left\vert l_{q}(F_{0})\right\vert =\left\Vert F_{0}\right\Vert .$
Hence, $\left\Vert l_{q}\right\Vert =1$. Recall that such a function $F_{0}$
exists on the basis of Urysohn's great lemma.

Note that if $q$ is a closed path, then $l_{q}$ annihilates all members of
the class $\mathcal{R}$. But in general, when $q$ is not closed, we do not
have the equality $l_{q}(G)=0,$ for all members $G\in \mathcal{R}$.
Nonetheless, this functional has the important property that

\begin{equation*}
\left\vert l_{q}(g_{1}+g_{2})\right\vert \leq \frac{2}{n}(\left\Vert
g_{1}\right\Vert +\left\Vert g_{2}\right\Vert ),\eqno(2.1)
\end{equation*}%
where $g_{1}$ and $g_{2}$ are ridge functions with the directions $\mathbf{a}%
_{1}$ and $\mathbf{a}_{2}$, respectively, that is, $g_{1}=g_{1}(\mathbf{a}%
_{1}\cdot \mathbf{x})$ and $g_{2}=g_{2}(\mathbf{a}_{2}\cdot \mathbf{x}).$
This property is important in the sense that if $n$ is sufficiently large,
then the functional $l_{q}$ is close to an annihilating functional. To prove
(2.1), note that $\left\vert l_{q}(g_{1})\right\vert \leq \frac{2}{n}%
\left\Vert g_{1}\right\Vert $ and $\left\vert l_{q}(g_{2})\right\vert \leq
\frac{2}{n}\left\Vert g_{2}\right\Vert $. These estimates become obvious if
consider the chain of equalities $g_{1}(\mathbf{a}_{1}\cdot \mathbf{x}%
_{1})=g_{1}(\mathbf{a}_{1}\cdot \mathbf{x}_{2}),$ $g_{1}(\mathbf{a}_{1}\cdot
\mathbf{x}_{3})=g_{1}(\mathbf{a}_{1}\cdot \mathbf{x}_{4}),...$(or $g_{1}(%
\mathbf{a}_{1}\cdot \mathbf{x}_{2})=g_{1}(\mathbf{a}_{1}\cdot \mathbf{x}%
_{3}),$ $g_{1}(\mathbf{a}_{1}\cdot \mathbf{x}_{4})=g_{1}(\mathbf{a}_{1}\cdot
\mathbf{x}_{5}),...$) for $g_{1}(\mathbf{a}_{1}\cdot \mathbf{x})$ and the
corresponding chain of equalities for $g_{2}(\mathbf{a}_{2}\cdot \mathbf{x})$%
%.

Now consider the infinite path $p=(\mathbf{p}_{1},\mathbf{p}_{2},...)$ and
form the finite paths $p_{k}=(\mathbf{p}_{1},...,\mathbf{p}_{k}),$ $%
k=1,2,... $. For ease of notation, let us set $l_{k}=l_{p_{k}}.$ The
sequence $\{l_{_{k}}\}_{k=1}^{\infty }$ is a subset of the unit ball of the
conjugate space $C^{\ast }(Q).$ By the Banach-Alaoglu theorem, the unit ball
is weak$^{\text{*}}$ compact in the weak$^{\text{*}}$ topology of $C^{\ast
}(Q)$ (see, e.g., Rudin \cite[p. 66]{Rud}). From this theorem we derive that
the sequence $\{l_{_{k}}\}_{k=1}^{\infty }$ must have weak$^{\text{*}}$
cluster points. Suppose $l^{\ast }$ denotes one of them. Without loss of
generality we may assume that $l_{k}\overset{weak^{\ast }}{\longrightarrow }%
l^{\ast },$ as $k\rightarrow \infty .$ From (2.1) it follows that $l^{\ast
}(g_{1}+g_{2})=0.$ That is, $l^{\ast }\in \mathcal{R}^{\bot },$ where the
symbol $\mathcal{R}^{\bot }$ stands for the annihilator of $\mathcal{R}$.
Since in addition $\left\Vert l^{\ast }\right\Vert \leq 1,$ we can write that

\begin{equation*}
\left\vert l^{\ast }(f)\right\vert =\left\vert l^{\ast }(f-G)\right\vert
\leq \left\Vert f-G\right\Vert ,\eqno(2.2)
\end{equation*}%
for all functions $G\in \mathcal{R}.$ On the other hand, since the infinite
bolt $p$ is extremal for $f-G_{0}$

\begin{equation*}
\left\vert l_{k}(f-G_{0})\right\vert =\left\Vert f-G_{0}\right\Vert ,\text{ }%
k=1,2,...
\end{equation*}%
Therefore,

\begin{equation*}
\left\vert l^{\ast }(f)\right\vert =\left\vert l^{\ast }(f-G_{0})\right\vert
=\left\Vert f-G_{0}\right\Vert .\eqno(2.3)
\end{equation*}%
From (2.2) and (2.3) we conclude that

\begin{equation*}
\left\Vert f-G_{0}\right\Vert \leq \left\Vert f-G\right\Vert ,
\end{equation*}%
for all $G\in \mathcal{R}.$ In other words, $G_{0}$ is a best approximation
to $f$. We proved the sufficiency of the theorem.

\textbf{Necessity.} The proof of this part is mainly based on the following
theorem of Singer.

\textbf{Theorem 2.2 (see Singer \cite{S}).} Let $X$ be a compact space, $U$
be a linear subspace of $C(X)$, $f\in C(X)\backslash U$ and $u_{0}\in U.$
Then $u_{0}$ is a best approximation to $f$ if and only if there exists a
regular Borel measure $\mu $ on $X$ such that

(1) The total variation $\left\Vert \mu \right\Vert =1$;

(2) $\mu $ is orthogonal to the subspace $U$, that is, $\int_{X}ud\mu =0$
for all $u\in U$;

(3) For the Jordan decomposition $\mu =\mu ^{+}-\mu ^{-}$,
\begin{equation*}
f(x)-u_{0}(x)=\left\{
\begin{array}{c}
\left\Vert f-u_{0}\right\Vert \text{ for }x\in S^{+}\text{,} \\
-\left\Vert f-u_{0}\right\Vert \text{ for }x\in S^{-}\text{,}%
\end{array}%
\right.
\end{equation*}%
where $S^{+}$ and $S^{-}$ are closed supports of the positive measures $\mu
^{+}$ and $\mu ^{-}$, respectively.

Let us show how we use this theorem in the proof of necessity part of our
theorem. Assume $G_{0}\in \mathcal{R}$ is a best approximation. For the
subspace $\mathcal{R},$ the existence of a measure $\mu $ satisfying the
conditions (1)-(3) is a direct consequence of Theorem 2.2. Let $\mathbf{x}%
_{0}$ be any point in $S^{+}.$\ Consider the point $y_{0}=\mathbf{a}%
_{1}\cdot \mathbf{x}_{0}$ and a $\delta $-neighborhood of $y_{0}$. That is,
choose an arbitrary $\delta >0$ and consider the set $I_{\delta
}=(y_{0}-\delta ,y_{0}+\delta )\cap \mathbf{a}_{1}\cdot Q.$ Here, $\mathbf{a}%
_{1}\cdot Q=\{\mathbf{a}_{1}\cdot \mathbf{x}:$ $\mathbf{x}\in Q\}.$ For any
subset $E\subset \mathbb{R}$, put

\begin{equation*}
E^{i}=\{\mathbf{x}\in Q:\mathbf{a}_{i}\cdot \mathbf{x}\in E\},\text{ }i=1,2.%
\text{ }
\end{equation*}

Clearly, for some sets $E,$ one or both the sets $E^{i}$ may be empty. Since
$I_{\delta }^{1}\cap S^{+}$ is not empty (note that $\mathbf{x}_{0}\in
I_{\delta }^{1}$), it follows that $\mu ^{+}(I_{\delta }^{1})>0.$ At the
same time $\mu (I_{\delta }^{1})=0,$ since $\mu $ is orthogonal to all
functions $g_{1}(\mathbf{a}_{1}\cdot \mathbf{x}).$ Therefore, $\mu
^{-}(I_{\delta }^{1})>0.$ We conclude that $I_{\delta }^{1}\cap S^{-}$ is
not empty. Denote this intersection by $A_{\delta }.$ Tending $\delta $ to $%
0,$ we obtain a set $A$ which is a subset of $S^{-}$ and has the property
that for each $\mathbf{x}\in A,$ we have $\mathbf{a}_{1}\cdot \mathbf{x}=%
\mathbf{a}_{1}\cdot \mathbf{x}_{0}.$ Fix any point $\mathbf{x}_{1}\in A$.
Changing $\mathbf{a}_{1}$, $\mu ^{+}$, $S^{+}$ to $\mathbf{a}_{2}$, $\mu ^{-}
$ and $S^{-}$ correspondingly, repeat the above process with the point $%
y_{1}=\mathbf{a}_{2}\cdot \mathbf{x}_{1}$ and a $\delta $-neighborhood of $%
y_{1}$. Then we obtain a point $\mathbf{x}_{2}\in S^{+}$ such that $\mathbf{a%
}_{2}\cdot \mathbf{x}_{2}=\mathbf{a}_{2}\cdot \mathbf{x}_{1}.$ Continuing
this process, one can construct points $\mathbf{x}_{3}$, $\mathbf{x}_{4}$,
and so on. Note that \ the set of all constructed points $\mathbf{x}_{i}$, $%
i=0,1,...,$ forms a path. By Theorem 2.2, this path is extremal for the
function $f-G_{0}$. We have proved the necessity and hence Theorem 2.1.

\bigskip

\textbf{Remark.} Theorem 2.1 was proven by Ismailov \cite{7} and in a more
general form by Pinkus \cite{P1} under additional assumption that $Q$ is
convex. Convexity assumption was made to guarantee continuity of the
following functions

\begin{equation*}
g_{1,i}(t)=\max_{\substack{ \mathbf{x}\in Q \\ \mathbf{a}_{i}\cdot \mathbf{x}%
=t}}F(\mathbf{x})\ \ \text{and }\ g_{2,i}(t)=\min\limits_{\substack{ \mathbf{%
x}\in Q \\ \mathbf{a}_{i}\cdot \mathbf{x}=t}}F(\mathbf{x}),\text{ }i=1,2,
\end{equation*}%
where $F$ is an arbitrary continuous function on $Q$. Note that in the proof
given above we need not continuity of these functions.

\bigskip

\end{document}